\def\bu{\bullet}
\def\marker{\>\hbox{${\vcenter{\vbox{
    \hrule height 0.4pt\hbox{\vrule width 0.4pt height 6pt
    \kern6pt\vrule width 0.4pt}\hrule height 0.4pt}}}$}\>}
\def\gpic#1{#1
     \medskip\par\noindent{\centerline{\box\graph}} \medskip}
\documentclass[12pt]{article}

\usepackage[margin=2cm]{geometry}

\usepackage{amsmath, amssymb, latexsym, amsthm, setspace}
\newtheorem{theorem}{Theorem} %Defines \begin{theorem} to write "Theorem"
\newtheorem{lemma}[theorem]{Lemma}

\newtheorem*{theorem*}{Theorem}
\newtheorem*{conjecture*}{Conjecture}
\newtheorem*{corollary*}{Corollary}
\newtheorem*{proposition*}{Proposition}

\theoremstyle{definition}

\newtheorem{definition}[theorem]{Definition}
\newtheorem{example}{Example}

\theoremstyle{remark}

\def\st{\colon\,}
\def\VEC#1#2#3{#1_{#2},\ldots,#1_{#3}}
\def\bT{{\bf T}}
\def\la{\langle}
\def\ra{\rangle}
\def\qed{\hfill$\Box$}

\usepackage[dvips]{graphicx}

\graphicspath{{counterfigs/}}

\author{
Armen S. Asratian\thanks{
Link\"oping University, Link\"oping Sweden, arasr@mai.liu.se.}\,,
Carl Johan Casselgren\thanks{
Ume\aa\enskip University, Ume\aa\enskip, Sweden, carl-johan.casselgren@math.umu.se.}\,,\\
Jennifer Vandenbussche\thanks{
University of Illinois, Urbana, IL, jarobin1@math.uiuc.edu.}\,,
Douglas B. West\thanks{
University of Illinois, Urbana, IL, west@math.uiuc.edu.  Work supported in part
by the NSA under Award No.~H98230-06-1-0065.}
}

\title{Proper path-factors and interval edge-coloring of $(3,4)$-biregular
bigraphs}
\date{April 6, 2007}

\begin{document}
\maketitle

\begin{abstract}
An {\it interval coloring} of a graph $G$ is a proper coloring of $E(G)$ by
positive integers such that the colors on the edges incident to any vertex are
consecutive.  A {\it $(3,4)$-biregular bigraph} is a bipartite graph in which
each vertex of one part has degree 3 and each vertex of the other has degree 4;
it is unknown whether these all have interval colorings.  We prove that $G$ has
an interval coloring using 6 colors when $G$ is a $(3,4)$-biregular bigraph
having a spanning subgraph whose components are paths with endpoints at
3-valent vertices and lengths in $\{2,4,6,8\}$.  We provide sufficient
conditions for the existence of such a subgraph.


AMSclass: 05C15, 05C70
\end{abstract}

\section{Introduction}

An {\it interval coloring} or {\it consecutive coloring} of a graph $G$ is a
proper coloring of the edges of $G$ by positive integers such that the colors
on the edges incident to any vertex are consecutive.  The notion was introduced
by Asratian and Kamalian~\cite{AK1} (available in English as~\cite{AK2}),
motivated by the problem of constructing timetables without ``gaps'' for
teachers and classes.  Hansen~\cite{Hansen} suggested another scenario:  a
school wishes to schedule parent-teacher conferences in time slots so that
every person's conferences occur in consecutive slots.  A solution
exists if and only if the bipartite graph with vertices for the people and
edges for the required meetings has an interval coloring.

In the context of edge-colorings, and particularly edge-colorings of bipartite
graphs, it is common to consider the general model in which multiple edges are
allowed.  In this paper, we adopt the convention that ``graph'' allows
multiple edges, and we will explicitly exclude multiple edges when necessary
(a {\it simple graph} is a graph without loops or multiple edges).

All regular bipartite graphs have interval colorings, since they have proper
edge-colorings in which all color classes are perfect matchings.
Not every graph has an interval coloring, since a graph $G$ with an interval
coloring must have a proper $\Delta(G)$-edge-coloring \cite{AK2}.  Furthermore,
Sevastjanov~\cite{Sevast} proved that determining whether a bipartite graph has
an interval coloring is NP-complete.  Nevertheless, trees
\cite{Hansen,Kamalian}, complete bipartite graphs~\cite{Hansen,Kamalian},
``doubly convex'' bipartite graphs~\cite{Kamalian}, grids~\cite{GK1}, and
simple outerplanar bipartite graphs~\cite{GK,Axe} all have interval colorings.
Giaro~\cite{Giaro} showed that one can decide in polynomial time whether
bipartite graphs with maximum degree 4 have interval 4-colorings.

An {\it $(a,b)$-biregular bigraph} is a bipartite graph where the vertices
in one part all have degree $a$ and the vertices in the other part all have
degree $b$.  Hansen~\cite{Hansen} proved that $(2,b)$-biregular bigraphs
are interval colorable when $b$ is even.  This was extended to all $b$ by
Hanson, Loten, and Toft~\cite{HansonLotenToft} and independently by
Kostochka~\cite{Kos}.  Kamalian~\cite{Kamalian} showed that the complete
bipartite graph $K_{b,a}$ has an interval coloring using $t$ colors if and only
if $a+b-\gcd(a,b)\le t\le a+b-1$, where $\gcd$ denotes the greatest common
divisor.  Asratian and Casselgren~\cite{AsCas} showed that recognizing whether 
$(3,6)$-biregular bigraphs have interval 6-colorings is NP-complete.

It is unknown whether all $(3,4)$-biregular bigraphs have interval colorings.
Hanson and Loten~\cite{HansonLoten} proved that no $(a,b)$-biregular bigraph
has an interval coloring with fewer than $a+b-\gcd(a,b)$ colors; thus
$(3,4)$-biregular bigraphs need at least 6 colors.  An {\it $X,Y$-bigraph} is a
bipartite graph with partite sets $X$ and $Y$.  In our $(3,4)$-biregular
$X,Y$-bigraphs, the vertices of $X$ will have degree 3.  Pyatkin~\cite{Pyatkin}
proved that if a $(3,4)$-biregular bigraph has a 3-regular subgraph covering the
vertices of degree 4, then it has an interval 6-coloring.

Here we obtain another sufficient condition for the existence of an interval
6-coloring of a $(3,4)$-biregular $X,Y$-bigraph $G$: If $G$ has a spanning
subgraph whose components are paths with endpoints in $X$ and lengths
in $\{2,4,6,8\}$ (we call this a {\it proper path-factor} of $G$), then
$G$ has an interval 6-coloring.  A longer proof of this was found earlier by
Casselgren~\cite{Cass}.

We present infinitely many $(3,4)$-biregular bigraphs that have proper
path-factors but do not satisfy Pyatkin's condition.  On the other hand,
$(3,4)$-biregular bigraphs with multiple edges need not have proper
path-factors, even if they satisfy Pyatkin's condition.  For example, consider
the graph formed from three triple-edges by adding a claw; that is,
the pairs $x_iy_i$ have multiplicity three for $i\in\{1,2,3\}$, and there is
an additional vertex $x_0$ with neighborhood $\{y_1,y_2,y_3\}$.  A 3-regular
subgraph covers $\{y_1,y_2,y_3\}$, but there is no proper path-factor.
Therefore, neither our result nor Pyatkin's result implies the other.

Various difficulties disappear when multiple edges are forbidden.  We have
found no simple $(3,4)$-biregular bigraph that does not have a proper
path-factor.  We conjecture that every simple $(3,4)$-biregular bigraph has a
proper path-factor.  In Section 3 we present various sufficient conditions for
the existence of a proper path-factor in such a graph.

\section{Interval 6-Colorings from Proper Path-Factors}

In general, an $\mathcal{H}$-factor of a graph is a spanning subgraph whose
components lie in $\mathcal{H}$.  We are interested in a particular family 
$\mathcal{H}$.  Let $d_H(v)$ denote the degree of a vertex $v$ in a graph $H$.

\begin{definition}
A {\it proper path-factor} of a $(3,4)$-biregular $X,Y$-bigraph $G$ is a
spanning subgraph of $G$ whose components are paths with endpoints in $X$
and lengths in $\{2,4,6,8\}$.
%A {\it $k$-valent} vertex is a vertex of degree $k$.
\end{definition}

Henceforth let $G$ be a $(3,4)$-biregular $X,Y$-bigraph.  Given a proper
path-factor $P$ of $G$, let $Q=G-E(P)$.  Observe that $d_Q(y)=2$ for all
$y\in Y$.  Furthermore, $d_Q(x)=2$ if $x$ is an endpoint of a component of $P$,
and $d_Q(x)=1$ if $x\in X$ and $x$ is an internal vertex of a component of $P$.
Thus every component of $Q$ is an even cycle or is a path with endpoints in
$X$.

\begin{definition}
Given a proper path-factor $P$ of $G$, the {\it $P$-graph of $G$}, denoted
$G_P$, is the graph with vertices $\{x \in X\st d_P(x)=2\}$ having
$x_i$ and $x_j$ adjacent when any condition below holds: \\
\null\quad (a) $x_i$ and $x_j$ are vertices of degree 2 in one component of $P$
with length 6, or \\
\null\quad (b) $x_i$ and $x_j$ are vertices of degree 2 at distance 4 in one
component of $P$ with length 8, or \\
\null\quad (c) $x_i$ and $x_j$ are vertices of degree 1 in one component of $Q$.
\end{definition}

\begin{lemma}
If $P$ is a proper path-factor of $G$, then $G_P$ is bipartite.  
\end{lemma}
\begin{proof}
Every vertex of $G_P$ has exactly one incident edge of type (c).  Some vertices
have one more neighbor, via an edge of type (a) or (b).  Thus
$\Delta(G_P)\le 2$.  Furthermore, the edges along any path or cycle in $G_P$
alternate type (c) with type (a) or (b).  Thus $G_P$ has no odd cycle.
\end{proof}

We say that a color appears ``at'' a vertex if it appears on an edge incident
to that vertex.

\begin{theorem}\label{pathfac}
If $G$ has a proper path-factor, then $G$ has an interval 6-coloring.
\end{theorem}
\begin{proof}
Let $P$ be a proper path-factor of $G$.  Let $c$ be a proper 2-coloring of
$V(G_P)$ with colors $A$ and $B$.  We define a 6-coloring of $E(G)$ that we
will show is an interval coloring.  Edges of $P$ receive colors from
$\{1,2,5,6\}$; edges of $Q$ receive colors from $\{3,4\}$.

First we color $E(Q)$.  Properly color cycles arbitrarily using colors 3 and 4.
A component of $Q$ that is a path has both endpoints in $G_P$, and they are
adjacent in $G_P$.  Hence $c(x)=A$ for one endpoint $x$ of the path, and
$c(x')=B$ for the other endpoint $x'$.  Alternate colors along the path,
starting with color 3 on the edge at $x$ and ending with color 4 on the edge at
$x'$.  Colors $3$ and $4$ both now appear at every vertex of $G$ having degree
2 in $Q$.

The edges of every component of $P$ are colored by alternating 2 and 1
(starting with 2) from one end, and alternating 5 and 6 (starting with 5)
from the other end.  We must specify which end is which and where to switch
from using one pair of colors to using the other.  The choice is based on the
colors that $c$ assigns to the internal vertices of the path that lie in $X$,
as illustrated in Figure 1.  Those vertices all have degree 1 in $Q$; they
appear in $V(G_P)$ and have colors under $c$.

\begin{figure}[hbt]
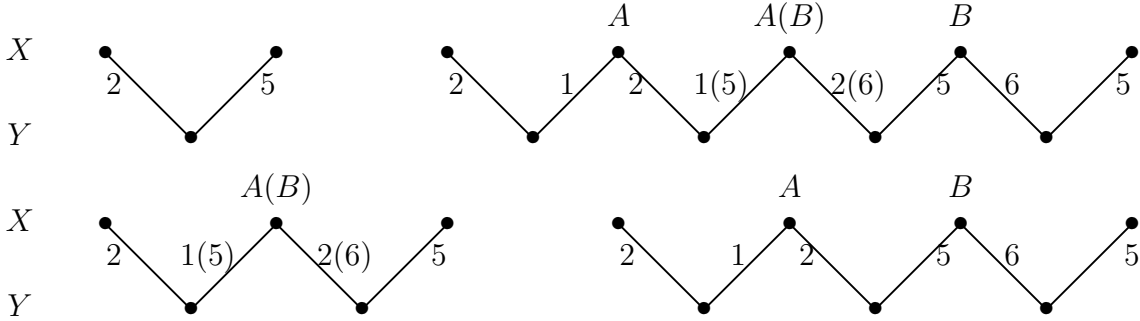

\gpic{
\expandafter\ifx\csname graph\endcsname\relax \csname newbox\endcsname\graph\fi
\expandafter\ifx\csname graphtemp\endcsname\relax \csname newdimen\endcsname\graphtemp\fi
\setbox\graph=\vtop{\vskip 0pt\hbox{%
    \graphtemp=.5ex\advance\graphtemp by 1.075in
    \rlap{\kern 0.000in\lower\graphtemp\hbox to 0pt{\hss $X$\hss}}%
    \graphtemp=.5ex\advance\graphtemp by 1.522in
    \rlap{\kern 0.000in\lower\graphtemp\hbox to 0pt{\hss $Y$\hss}}%
    \graphtemp=.5ex\advance\graphtemp by 1.075in
    \rlap{\kern 0.448in\lower\graphtemp\hbox to 0pt{\hss $\bu$\hss}}%
    \graphtemp=.5ex\advance\graphtemp by 1.522in
    \rlap{\kern 0.896in\lower\graphtemp\hbox to 0pt{\hss $\bu$\hss}}%
    \graphtemp=.5ex\advance\graphtemp by 1.075in
    \rlap{\kern 1.343in\lower\graphtemp\hbox to 0pt{\hss $\bu$\hss}}%
    \graphtemp=.5ex\advance\graphtemp by 1.522in
    \rlap{\kern 1.791in\lower\graphtemp\hbox to 0pt{\hss $\bu$\hss}}%
    \graphtemp=.5ex\advance\graphtemp by 1.075in
    \rlap{\kern 2.239in\lower\graphtemp\hbox to 0pt{\hss $\bu$\hss}}%
    \special{pn 11}%
    \special{pa 448 1075}%
    \special{pa 896 1522}%
    \special{fp}%
    \special{pa 896 1522}%
    \special{pa 1343 1075}%
    \special{fp}%
    \special{pa 1343 1075}%
    \special{pa 1791 1522}%
    \special{fp}%
    \special{pa 1791 1522}%
    \special{pa 2239 1075}%
    \special{fp}%
    \graphtemp=.5ex\advance\graphtemp by 1.254in
    \rlap{\kern 0.493in\lower\graphtemp\hbox to 0pt{\hss 2\hss}}%
    \graphtemp=.5ex\advance\graphtemp by 1.254in
    \rlap{\kern 0.985in\lower\graphtemp\hbox to 0pt{\hss 1(5)\hss}}%
    \graphtemp=.5ex\advance\graphtemp by 1.254in
    \rlap{\kern 1.701in\lower\graphtemp\hbox to 0pt{\hss 2(6)\hss}}%
    \graphtemp=.5ex\advance\graphtemp by 1.254in
    \rlap{\kern 2.194in\lower\graphtemp\hbox to 0pt{\hss 5\hss}}%
    \graphtemp=.5ex\advance\graphtemp by 0.896in
    \rlap{\kern 1.343in\lower\graphtemp\hbox to 0pt{\hss $A(B)$\hss}}%
    \graphtemp=.5ex\advance\graphtemp by 1.075in
    \rlap{\kern 3.134in\lower\graphtemp\hbox to 0pt{\hss $\bu$\hss}}%
    \graphtemp=.5ex\advance\graphtemp by 1.522in
    \rlap{\kern 3.582in\lower\graphtemp\hbox to 0pt{\hss $\bu$\hss}}%
    \graphtemp=.5ex\advance\graphtemp by 1.075in
    \rlap{\kern 4.030in\lower\graphtemp\hbox to 0pt{\hss $\bu$\hss}}%
    \graphtemp=.5ex\advance\graphtemp by 1.522in
    \rlap{\kern 4.478in\lower\graphtemp\hbox to 0pt{\hss $\bu$\hss}}%
    \graphtemp=.5ex\advance\graphtemp by 1.075in
    \rlap{\kern 4.925in\lower\graphtemp\hbox to 0pt{\hss $\bu$\hss}}%
    \graphtemp=.5ex\advance\graphtemp by 1.522in
    \rlap{\kern 5.373in\lower\graphtemp\hbox to 0pt{\hss $\bu$\hss}}%
    \graphtemp=.5ex\advance\graphtemp by 1.075in
    \rlap{\kern 5.821in\lower\graphtemp\hbox to 0pt{\hss $\bu$\hss}}%
    \special{pa 3134 1075}%
    \special{pa 3582 1522}%
    \special{fp}%
    \special{pa 3582 1522}%
    \special{pa 4030 1075}%
    \special{fp}%
    \special{pa 4030 1075}%
    \special{pa 4478 1522}%
    \special{fp}%
    \special{pa 4478 1522}%
    \special{pa 4925 1075}%
    \special{fp}%
    \special{pa 4925 1075}%
    \special{pa 5373 1522}%
    \special{fp}%
    \special{pa 5373 1522}%
    \special{pa 5821 1075}%
    \special{fp}%
    \graphtemp=.5ex\advance\graphtemp by 1.254in
    \rlap{\kern 3.179in\lower\graphtemp\hbox to 0pt{\hss 2\hss}}%
    \graphtemp=.5ex\advance\graphtemp by 1.254in
    \rlap{\kern 3.761in\lower\graphtemp\hbox to 0pt{\hss 1\hss}}%
    \graphtemp=.5ex\advance\graphtemp by 1.254in
    \rlap{\kern 4.119in\lower\graphtemp\hbox to 0pt{\hss 2\hss}}%
    \graphtemp=.5ex\advance\graphtemp by 1.254in
    \rlap{\kern 4.836in\lower\graphtemp\hbox to 0pt{\hss 5\hss}}%
    \graphtemp=.5ex\advance\graphtemp by 1.254in
    \rlap{\kern 5.194in\lower\graphtemp\hbox to 0pt{\hss 6\hss}}%
    \graphtemp=.5ex\advance\graphtemp by 1.254in
    \rlap{\kern 5.821in\lower\graphtemp\hbox to 0pt{\hss 5\hss}}%
    \graphtemp=.5ex\advance\graphtemp by 0.896in
    \rlap{\kern 4.030in\lower\graphtemp\hbox to 0pt{\hss $A$\hss}}%
    \graphtemp=.5ex\advance\graphtemp by 0.896in
    \rlap{\kern 4.925in\lower\graphtemp\hbox to 0pt{\hss $B$\hss}}%
    \graphtemp=.5ex\advance\graphtemp by 0.179in
    \rlap{\kern 0.000in\lower\graphtemp\hbox to 0pt{\hss $X$\hss}}%
    \graphtemp=.5ex\advance\graphtemp by 0.627in
    \rlap{\kern 0.000in\lower\graphtemp\hbox to 0pt{\hss $Y$\hss}}%
    \graphtemp=.5ex\advance\graphtemp by 0.179in
    \rlap{\kern 0.448in\lower\graphtemp\hbox to 0pt{\hss $\bu$\hss}}%
    \graphtemp=.5ex\advance\graphtemp by 0.627in
    \rlap{\kern 0.896in\lower\graphtemp\hbox to 0pt{\hss $\bu$\hss}}%
    \graphtemp=.5ex\advance\graphtemp by 0.179in
    \rlap{\kern 1.343in\lower\graphtemp\hbox to 0pt{\hss $\bu$\hss}}%
    \special{pa 448 179}%
    \special{pa 896 627}%
    \special{fp}%
    \special{pa 896 627}%
    \special{pa 1343 179}%
    \special{fp}%
    \graphtemp=.5ex\advance\graphtemp by 0.358in
    \rlap{\kern 0.493in\lower\graphtemp\hbox to 0pt{\hss 2\hss}}%
    \graphtemp=.5ex\advance\graphtemp by 0.358in
    \rlap{\kern 1.299in\lower\graphtemp\hbox to 0pt{\hss 5\hss}}%
    \graphtemp=.5ex\advance\graphtemp by 0.179in
    \rlap{\kern 2.239in\lower\graphtemp\hbox to 0pt{\hss $\bu$\hss}}%
    \graphtemp=.5ex\advance\graphtemp by 0.627in
    \rlap{\kern 2.687in\lower\graphtemp\hbox to 0pt{\hss $\bu$\hss}}%
    \graphtemp=.5ex\advance\graphtemp by 0.179in
    \rlap{\kern 3.134in\lower\graphtemp\hbox to 0pt{\hss $\bu$\hss}}%
    \graphtemp=.5ex\advance\graphtemp by 0.627in
    \rlap{\kern 3.582in\lower\graphtemp\hbox to 0pt{\hss $\bu$\hss}}%
    \graphtemp=.5ex\advance\graphtemp by 0.179in
    \rlap{\kern 4.030in\lower\graphtemp\hbox to 0pt{\hss $\bu$\hss}}%
    \graphtemp=.5ex\advance\graphtemp by 0.627in
    \rlap{\kern 4.478in\lower\graphtemp\hbox to 0pt{\hss $\bu$\hss}}%
    \graphtemp=.5ex\advance\graphtemp by 0.179in
    \rlap{\kern 4.925in\lower\graphtemp\hbox to 0pt{\hss $\bu$\hss}}%
    \graphtemp=.5ex\advance\graphtemp by 0.627in
    \rlap{\kern 5.373in\lower\graphtemp\hbox to 0pt{\hss $\bu$\hss}}%
    \graphtemp=.5ex\advance\graphtemp by 0.179in
    \rlap{\kern 5.821in\lower\graphtemp\hbox to 0pt{\hss $\bu$\hss}}%
    \special{pa 2239 179}%
    \special{pa 2687 627}%
    \special{fp}%
    \special{pa 2687 627}%
    \special{pa 3134 179}%
    \special{fp}%
    \special{pa 3134 179}%
    \special{pa 3582 627}%
    \special{fp}%
    \special{pa 3582 627}%
    \special{pa 4030 179}%
    \special{fp}%
    \special{pa 4030 179}%
    \special{pa 4478 627}%
    \special{fp}%
    \special{pa 4478 627}%
    \special{pa 4925 179}%
    \special{fp}%
    \special{pa 4925 179}%
    \special{pa 5373 627}%
    \special{fp}%
    \special{pa 5373 627}%
    \special{pa 5821 179}%
    \special{fp}%
    \graphtemp=.5ex\advance\graphtemp by 0.358in
    \rlap{\kern 2.284in\lower\graphtemp\hbox to 0pt{\hss 2\hss}}%
    \graphtemp=.5ex\advance\graphtemp by 0.358in
    \rlap{\kern 2.866in\lower\graphtemp\hbox to 0pt{\hss 1\hss}}%
    \graphtemp=.5ex\advance\graphtemp by 0.358in
    \rlap{\kern 3.224in\lower\graphtemp\hbox to 0pt{\hss 2\hss}}%
    \graphtemp=.5ex\advance\graphtemp by 0.358in
    \rlap{\kern 5.776in\lower\graphtemp\hbox to 0pt{\hss 5\hss}}%
    \graphtemp=.5ex\advance\graphtemp by 0.358in
    \rlap{\kern 5.194in\lower\graphtemp\hbox to 0pt{\hss 6\hss}}%
    \graphtemp=.5ex\advance\graphtemp by 0.358in
    \rlap{\kern 4.836in\lower\graphtemp\hbox to 0pt{\hss 5\hss}}%
    \graphtemp=.5ex\advance\graphtemp by 0.000in
    \rlap{\kern 3.134in\lower\graphtemp\hbox to 0pt{\hss $A$\hss}}%
    \graphtemp=.5ex\advance\graphtemp by 0.000in
    \rlap{\kern 4.030in\lower\graphtemp\hbox to 0pt{\hss $A(B)$\hss}}%
    \graphtemp=.5ex\advance\graphtemp by 0.000in
    \rlap{\kern 4.925in\lower\graphtemp\hbox to 0pt{\hss $B$\hss}}%
    \graphtemp=.5ex\advance\graphtemp by 0.358in
    \rlap{\kern 3.672in\lower\graphtemp\hbox to 0pt{\hss 1(5)\hss}}%
    \graphtemp=.5ex\advance\graphtemp by 0.358in
    \rlap{\kern 4.388in\lower\graphtemp\hbox to 0pt{\hss 2(6)\hss}}%
    \hbox{\vrule depth1.701in width0pt height 0pt}%
    \kern 6.000in
  }%
}%
}
\caption{Coloring the edges of $P$}
\end{figure}

Let $H$ be a component of $P$.  If $H\cong P_3$, then we assign 2 and 5 to the
edges arbitrarily.
If $H\cong P_5$ with middle vertex $x$, then it does not matter which end
edge gets color 2 and which gets color 5, but the middle edges get colors 1 and
2 if $c(x)=A$, 5 and 6 if $c(x)=B$.
If $H\cong P_7$, then the internal vertices are adjacent in $G_P$ and receive
distinct colors under $c$; use $2,1,2$ from the end closest to the one colored
$A$ and $5,6,5$ from the end closest to the one colored $B$.
If $H\cong P_9$, then the internal vertices at distance 4 on the path again are
adjacent in $G_P$, and the three edges from each end are colored in the same
way as for $P_7$.  The two central edges are colored like the middle edges
of $P_5$, based on the color under $c$ of the central vertex of the path.

We check that the resulting $6$-edge-coloring is an interval coloring.
Each vertex of $Y$ has colors 3 and 4 on its incident edges in $Q$ and
receives $\{2,5\}$ or $\{1,2\}$ or $\{5,6\}$ on its incident edges in $P$,
forming an interval in each case.  Each endpoint of a component of $P$
has colors 3 and 4 from $Q$ and receives color 2 or 5 from $P$.
Each internal vertex $x$ of a component of $P$ receives 3 from $Q$ and   
$\{1,2\}$ from $P$ if $c(x)=A$, while it receives 4 from $Q$ and $\{5,6\}$
from $P$ if $c(x)=B$.
\end{proof}

\bigskip

This technique does not extend to arbitrary path and cycle factors.
We switch from $1,2$-alternation to $5,6$-alternation only once along a path in
$P$ and cannot switch back.  Thus we need that along any path of $P$, the
internal vertices with color $A$ under $c$ all precede those with color $B$.
With longer paths, our technique offers no mechanism for achieving this; the
graph $G_P$ can only enforce that vertices receive different colors under $c$.
Introducing more edges into $G_P$ to prevent alternation of $A$ and $B$
along the path destroys the 2-colorability of $G_P$.

\section{Constructions and Conditions for Proper Path-Factors}

To apply the theorem, we seek proper path-factors of $(3,4)$-biregular bigraphs.
Here we will give some sufficient conditions for existence of proper
path-factors and provide some examples related to Pyatkin's condition.

We call a 3-regular subgraph of a $(3,4)$-biregular bigraph that
covers the vertices of degree 4 a {\it full 3-regular subgraph}.  Pyatkin
proved that a $(3,4)$-biregular bigraph with a full 3-regular
subgraph has an interval 6-coloring.  We begin with an example that satisfies
our condition but not Pyatkin's condition.  Let $[n]=\{1,\ldots,n\}$.

\begin{example}\label{sets}
{\it The $X,Y$-bigraph $G$ defined by letting $X$ and $Y$ be the 3-sets and
2-sets in $[6]$, with adjacency defined by proper containment, has an interval
6-coloring.} By Theorem~\ref{pathfac}, it suffices to find a proper path-factor.
In fact, $G$ has a $P_7$-factor as shown below.

\begin{centering}

$124 \rightarrow 12 \rightarrow 123 \rightarrow 23 \rightarrow 235 \rightarrow 35 \rightarrow 345$

$135 \rightarrow 13 \rightarrow 134 \rightarrow 34 \rightarrow 346 \rightarrow 46 \rightarrow 456$

$146 \rightarrow 14 \rightarrow 145 \rightarrow 45 \rightarrow 245 \rightarrow 25 \rightarrow 256$

$125 \rightarrow 15 \rightarrow 156 \rightarrow 56 \rightarrow 356 \rightarrow 36 \rightarrow 236$

$136 \rightarrow 16 \rightarrow 126 \rightarrow 26 \rightarrow 246 \rightarrow 24 \rightarrow 234$

\end{centering}

Here $|X|=20$ and $|Y|=15$, with $Y$ corresponding to the edge set of $K_6$.
The neighborhood of a vertex in $X$ corresponds to a triangle in $K_6$.
Hence five vertices can be deleted from $G$ to leave a full 3-regular subgraph
if and only if $K_6$ decomposes into five triangles.  It does not, because the
vertices of $K_6$ have odd degree.
\qed
\end{example}

We next construct infinitely many examples that satisfy our condition but not
Pyatkin's, starting with a graph smaller than that of Example~\ref{sets}.

\begin{example}
The smallest simple $(3,4)$-biregular bigraph is $K_{3,4}$; it satisfies
Pyatkin's condition.  The next smallest such graphs have eight vertices of
degree 3 and six of degree 4.  For example, consider an $X,Y$-bigraph where
$Y=[6]$ and the neighborhoods of the vertices in $X$ are eight triples from
$[6]$, with each element used in four triples.  The graph fails Pyatkin's
condition if and only if the triple system does not have two disjoint triples.

Case analysis shows that it is not possible to avoid two disjoint triples
without a repeated triple.  However, it is possible using a repeated triple,
as in $\{123,124,235,346,346,145,156,256\}$.  The resulting $(3,4)$-biregular
bigraph has a $P_7$-factor as shown in bold in Figure 2.
\qed
\end{example}

\begin{figure}[bht]
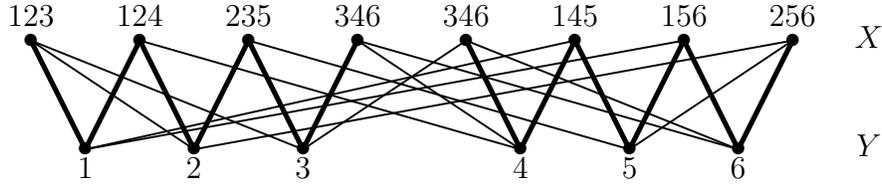

\gpic{
\expandafter\ifx\csname graph\endcsname\relax \csname newbox\endcsname\graph\fi
\expandafter\ifx\csname graphtemp\endcsname\relax \csname newdimen\endcsname\graphtemp\fi
\setbox\graph=\vtop{\vskip 0pt\hbox{%
    \graphtemp=.5ex\advance\graphtemp by 0.114in
    \rlap{\kern 0.114in\lower\graphtemp\hbox to 0pt{\hss $\bu$\hss}}%
    \graphtemp=.5ex\advance\graphtemp by 0.684in
    \rlap{\kern 0.399in\lower\graphtemp\hbox to 0pt{\hss $\bu$\hss}}%
    \graphtemp=.5ex\advance\graphtemp by 0.114in
    \rlap{\kern 0.684in\lower\graphtemp\hbox to 0pt{\hss $\bu$\hss}}%
    \graphtemp=.5ex\advance\graphtemp by 0.684in
    \rlap{\kern 0.968in\lower\graphtemp\hbox to 0pt{\hss $\bu$\hss}}%
    \graphtemp=.5ex\advance\graphtemp by 0.114in
    \rlap{\kern 1.253in\lower\graphtemp\hbox to 0pt{\hss $\bu$\hss}}%
    \graphtemp=.5ex\advance\graphtemp by 0.684in
    \rlap{\kern 1.538in\lower\graphtemp\hbox to 0pt{\hss $\bu$\hss}}%
    \graphtemp=.5ex\advance\graphtemp by 0.114in
    \rlap{\kern 1.823in\lower\graphtemp\hbox to 0pt{\hss $\bu$\hss}}%
    \graphtemp=.5ex\advance\graphtemp by 0.797in
    \rlap{\kern 0.399in\lower\graphtemp\hbox to 0pt{\hss 1\hss}}%
    \graphtemp=.5ex\advance\graphtemp by 0.797in
    \rlap{\kern 0.968in\lower\graphtemp\hbox to 0pt{\hss 2\hss}}%
    \graphtemp=.5ex\advance\graphtemp by 0.797in
    \rlap{\kern 1.538in\lower\graphtemp\hbox to 0pt{\hss 3\hss}}%
    \graphtemp=.5ex\advance\graphtemp by 0.000in
    \rlap{\kern 0.114in\lower\graphtemp\hbox to 0pt{\hss 123\hss}}%
    \graphtemp=.5ex\advance\graphtemp by 0.000in
    \rlap{\kern 0.684in\lower\graphtemp\hbox to 0pt{\hss 124\hss}}%
    \graphtemp=.5ex\advance\graphtemp by 0.000in
    \rlap{\kern 1.253in\lower\graphtemp\hbox to 0pt{\hss 235\hss}}%
    \graphtemp=.5ex\advance\graphtemp by 0.000in
    \rlap{\kern 1.823in\lower\graphtemp\hbox to 0pt{\hss 346\hss}}%
    \special{pn 28}%
    \special{pa 114 114}%
    \special{pa 399 684}%
    \special{fp}%
    \special{pa 399 684}%
    \special{pa 684 114}%
    \special{fp}%
    \special{pa 684 114}%
    \special{pa 968 684}%
    \special{fp}%
    \special{pa 968 684}%
    \special{pa 1253 114}%
    \special{fp}%
    \special{pa 1253 114}%
    \special{pa 1538 684}%
    \special{fp}%
    \special{pa 1538 684}%
    \special{pa 1823 114}%
    \special{fp}%
    \graphtemp=.5ex\advance\graphtemp by 0.114in
    \rlap{\kern 2.392in\lower\graphtemp\hbox to 0pt{\hss $\bu$\hss}}%
    \graphtemp=.5ex\advance\graphtemp by 0.684in
    \rlap{\kern 2.677in\lower\graphtemp\hbox to 0pt{\hss $\bu$\hss}}%
    \graphtemp=.5ex\advance\graphtemp by 0.114in
    \rlap{\kern 2.962in\lower\graphtemp\hbox to 0pt{\hss $\bu$\hss}}%
    \graphtemp=.5ex\advance\graphtemp by 0.684in
    \rlap{\kern 3.247in\lower\graphtemp\hbox to 0pt{\hss $\bu$\hss}}%
    \graphtemp=.5ex\advance\graphtemp by 0.114in
    \rlap{\kern 3.532in\lower\graphtemp\hbox to 0pt{\hss $\bu$\hss}}%
    \graphtemp=.5ex\advance\graphtemp by 0.684in
    \rlap{\kern 3.816in\lower\graphtemp\hbox to 0pt{\hss $\bu$\hss}}%
    \graphtemp=.5ex\advance\graphtemp by 0.114in
    \rlap{\kern 4.101in\lower\graphtemp\hbox to 0pt{\hss $\bu$\hss}}%
    \graphtemp=.5ex\advance\graphtemp by 0.797in
    \rlap{\kern 2.677in\lower\graphtemp\hbox to 0pt{\hss 4\hss}}%
    \graphtemp=.5ex\advance\graphtemp by 0.797in
    \rlap{\kern 3.247in\lower\graphtemp\hbox to 0pt{\hss 5\hss}}%
    \graphtemp=.5ex\advance\graphtemp by 0.797in
    \rlap{\kern 3.816in\lower\graphtemp\hbox to 0pt{\hss 6\hss}}%
    \graphtemp=.5ex\advance\graphtemp by 0.000in
    \rlap{\kern 2.392in\lower\graphtemp\hbox to 0pt{\hss 346\hss}}%
    \graphtemp=.5ex\advance\graphtemp by 0.000in
    \rlap{\kern 2.962in\lower\graphtemp\hbox to 0pt{\hss 145\hss}}%
    \graphtemp=.5ex\advance\graphtemp by 0.000in
    \rlap{\kern 3.532in\lower\graphtemp\hbox to 0pt{\hss 156\hss}}%
    \graphtemp=.5ex\advance\graphtemp by 0.000in
    \rlap{\kern 4.101in\lower\graphtemp\hbox to 0pt{\hss 256\hss}}%
    \special{pa 2392 114}%
    \special{pa 2677 684}%
    \special{fp}%
    \special{pa 2677 684}%
    \special{pa 2962 114}%
    \special{fp}%
    \special{pa 2962 114}%
    \special{pa 3247 684}%
    \special{fp}%
    \special{pa 3247 684}%
    \special{pa 3532 114}%
    \special{fp}%
    \special{pa 3532 114}%
    \special{pa 3816 684}%
    \special{fp}%
    \special{pa 3816 684}%
    \special{pa 4101 114}%
    \special{fp}%
    \graphtemp=.5ex\advance\graphtemp by 0.114in
    \rlap{\kern 4.500in\lower\graphtemp\hbox to 0pt{\hss $X$\hss}}%
    \graphtemp=.5ex\advance\graphtemp by 0.684in
    \rlap{\kern 4.500in\lower\graphtemp\hbox to 0pt{\hss $Y$\hss}}%
    \special{pn 11}%
    \special{pa 2962 114}%
    \special{pa 399 684}%
    \special{fp}%
    \special{pa 399 684}%
    \special{pa 3532 114}%
    \special{fp}%
    \special{pa 114 114}%
    \special{pa 968 684}%
    \special{fp}%
    \special{pa 968 684}%
    \special{pa 4101 114}%
    \special{fp}%
    \special{pa 114 114}%
    \special{pa 1538 684}%
    \special{fp}%
    \special{pa 1538 684}%
    \special{pa 2392 114}%
    \special{fp}%
    \special{pa 684 114}%
    \special{pa 2677 684}%
    \special{fp}%
    \special{pa 2677 684}%
    \special{pa 1823 114}%
    \special{fp}%
    \special{pa 1253 114}%
    \special{pa 3247 684}%
    \special{fp}%
    \special{pa 3247 684}%
    \special{pa 4101 114}%
    \special{fp}%
    \special{pa 1823 114}%
    \special{pa 3816 684}%
    \special{fp}%
    \special{pa 3816 684}%
    \special{pa 2392 114}%
    \special{fp}%
    \hbox{\vrule depth0.797in width0pt height 0pt}%
    \kern 4.500in
  }%
}%
}
\caption{$P_7$-factor in a bigraph with no full 3-regular subgraph}
\end{figure}

%\eject

Using the next lemma, we can generate infinitely many examples that have
$P_7$-factors but have no full 3-regular subgraphs.  The number of vertices can
be any nontrivial multiple of $7$.  Here multiple edges are allowed.

\begin{lemma}\label{2switch}
For $i\in\{1,2\}$, let $G_i$ be a 2-edge-connected $(3,4)$-biregular bigraph
having a $P_7$-factor $F_i$, and choose $e_i\in E(G_i)-E(F_i)$.  Let
$G$ be the $(3,4)$-biregular bigraph obtained from the disjoint
union of $G_1$ and $G_2$ by deleting $e_1$ and $e_2$ and replacing them with
two other edges $e_1'$ and $e_2'$ joining their endpoints.  If $G_1$ has no
full 3-regular subgraph, then $G$ is a larger 2-edge-connected
$(3,4)$-biregular bigraph having a $P_7$-factor but no full
3-regular subgraph.
\end{lemma}
\begin{proof}
Since $e_i\notin E(F_i)$, the subgraph $F_1\cup F_2$ is a $P_7$-factor of $G$.
Since each $G_i$ is 2-edge-connected, $G$ is connected.  Also, a cycle through
$e_i$ in $G_i$ can detour through $G_{3-i}$ using $e_1'$ and $e_2'$.  Thus
$G$ is 2-edge-connected.

Suppose that $G$ has a full 3-regular subgraph $H$.  By considering vertex
degrees, $H$ must have an even number of edges in $\{e_1',e_2'\}$.
If $H$ uses neither, then $H$ restricts to full 3-regular subgraphs of $G_1$
and $G_2$.  If $H$ uses both, then replacing $e_1'$ and $e_2'$ with $e_1$ and
$e_2$ yields full 3-regular subgraphs of $G_1$ and $G_2$.
\end{proof}

If $G_1$ and $G_2$ in Lemma~\ref{2switch} have no multiple edges, then neither
does the resulting graph $G$.

Our next theorem gives a sufficient condition for existence of a proper
path-factor in a $(3,4)$-biregular bigraph.  First we note an easy
lemma.

\begin{lemma}\label{p3fac}
Every $(2,4)$-biregular bigraph $H$ has a $(1,2)$-biregular factor
with every component isomorphic to $P_3$.  (Indeed, $H$ decomposes into
two such factors.)
\end{lemma}
\begin{proof}
Each component of $H$ is Eulerian and has an even number of edges.  Taking the
even-indexed edges from an Eulerian circuit in each component takes half the
edges from each vertex.  Thus it yields a spanning subgraph in which every
vertex of one partite set has degree 1 and every vertex of the other has degree
2.  Hence each component of the subgraph is isomorphic to $P_3$.
\end{proof}

\begin{theorem}
A $(3,4)$-biregular bigraph $G$ has a $P_7$-factor (and hence an
interval 6-coloring) if $G$ has a $(2,4)$-biregular subgraph covering the set
of vertices of degree 3.
\end{theorem}
\begin{proof}
Let $G$ have bipartition $X,Y$, where $|X| = 4k$ and $|Y| = 3k$.  Let $H$ be a
$(2,4)$-biregular subgraph of $G$ covering $X$; we obtain $H$ from $G$ by
deleting vertices $\VEC u1k$ of $Y$ that have disjoint neighborhoods.
Let $\hat Y=\{\VEC u1k\}$ and $Y'=Y-\hat Y$, so $H$ has bipartition $X,Y'$.

By Lemma~\ref{p3fac}, $H$ has a spanning subgraph $F$ whose components are
copies of $P_3$ with endpoints in $X$.  Let $\VEC T1{2k}$ be these paths.
Index $X$ so that $V(T_i)=\{x_{2i-1},y_i,x_{2i}\}$ (we maintain the flexibility
to decide later which end is $x_{2i-1}$ and which is $x_{2i}$).

Next we obtain from $G-Y'$ a graph $H'$ by combining the endpoints of each
path $T_i$ into a single vertex $x'_i$.  Since $G-Y'$ is a $(1,4)$-biregular
$X,\hat Y$-bigraph, $H'$ is a $(2,4)$-biregular $X',\hat Y$-bigraph, where
$X'=\{\VEC{x'}1{2k}\}$.  Note that multiple edges may arise in $H'$.

For each of the $k$ vertices of $\hat Y$, we construct a path of length $6$ in
$G$ with endpoints in $X$.  By Lemma~\ref{p3fac}, $H'$ has a spanning subgraph
$F'$ whose components are copies of $P_3$ with endpoints in $X'$.  For
$u\in\hat Y$, let $x'_i$ and $x'_j$ be the neighbors of $u$ in $F'$.  Thus in
$G$ the vertex $u$ is adjacent to one endpoint of $T_i$ and one endpoint of
$T_j$.  We may complete the indexing of $X$ so that these neighbors of $u$ are
$x_{2i}\in V(T_{i})$ and $x_{2j-1} \in V(T_{j})$.  The path we associate with
$u$ is then $\la x_{2j-1},y_{i},x_{2i},u,x_{2j-1},y_{j},x_{2j}\ra$, isomorphic
to $P_7$.

We check that these paths are pairwise disjoint.  Each uses exactly one vertex
of $\hat Y$.  Since $F'$ has exactly one edge incident to each vertex of $X'$,
for each $i$ the vertices of $T_i$ occur in exactly one of the paths.  Hence
these paths form a $P_7$-factor, and Theorem~\ref{pathfac} applies.
\end{proof}

We now return to simple $(3,4)$-biregular bigraphs.  Although the examples
constructed so far in this section all have $P_7$-factors,
Casselgren~\cite{Cass} found a simple $(3,4)$-biregular bigraph with no
$P_7$-factor.  We conjecture that every simple $(3,4)$-biregular bigraph
has the weaker property of having a proper path-factor.  It should also hold
that Pyatkin's condition guarantees the existence of a proper path-factor
in a simple $(3,4)$-biregular bigraph, but this also seems difficult.  We
present a condition that guarantees a proper path-factor when combined with
Pyatkin's condition.

Let $G$ be a simple $(3,4)$-biregular $X,Y$-bigraph having a full 3-regular
subgraph $H$.  Since $|X|=4k$ and $|Y|=3k$ for some $k$, we may let
$X'=X\cap V(H)$ and $X_0=X-X'$, where $X_0=\{\VEC{x^0}1k\}$.  Since $H$ is
3-regular, $H$ has a proper 3-edge coloring.  Fix such a coloring $c$, and let
$H'$ be the spanning subgraph of $H$ whose edges are those with color 1 or 2
under $c$.  Define an auxiliary graph $F$ with vertex set $Y$ by putting
$y_iy_j \in E(F)$ if $H'$ has a $y_i,y_j$-path of length~2.  Note that $F$ may
have multiple edges and is 2-regular, since each vertex of $Y$ has one incident
edge with each color under $c$.  Since $G$ is simple, the components of $F$ are
cycles of length at least~2.

Since $G$ is $(3,4)$-biregular, the neighborhoods of the vertices of $X_0$
partition $Y$ into triples; let $T_i=N_G(x^0_i)=\{y_i^1,y_i^2,y_i^3\}$.
Let $\bT$ denote the family $\VEC T1k$.

\begin{definition}\label{transv}
For families of disjoint triples, we define a {\it transversal} to be a set $S$
having exactly one element from each triple.  For a family $\bT$ defined on the
vertices of a 2-regular graph $F$, an {\it independent transversal} is a
transversal $S$ that is an independent set in $F$.  A {\it spread transversal}
is a transversal $S$ such that, given directions on the cycles of $F$, for
every vertex $v$ of $F$ that does not belong to $S$, there is a vertex of $S$
among the next three vertices after $v$ along the forward direction of its
cycle in $F$.  Let $F^*$ be the 4-regular graph obtained from $F$ by
adding triangles whose vertex sets are the triples of $\bT$.  A {\it mixed
transversal} is a transversal that restricts on each component of $F^*$ to an
independent transversal or a spread transversal.
\end{definition}

Note that a spread transversal intersects each cycle of $F$.

\begin{theorem}
Let $G$ be a simple $(3,4)$-biregular $X,Y$-bigraph having a full 3-regular
subgraph $H$, and let $F$ and $\bT$ be the $2$-regular graph and triple system
defined as above.  If $\bT$ has a mixed transversal, then $G$ has a proper
path-factor.
\end{theorem}

\begin{proof}
Let $c$ be a proper 3-edge-coloring of $H$, and let $M$ be the perfect matching
in color 1.  The $i$th triple in $\bT$ is $\{y_i^1,y_i^2,y_i^3\}$; we may let
$y_i^1$ be the vertex of $T_i$ in the mixed transversal $Y_1$.  Let
$Y_2=\{\VEC{y^2}1k\}$ and $Y_3=\{\VEC{y^3}1k\}$.  For $x\in X'=X-X_0$, put
$x\in X_j$ if the other endpoint of the edge of $M$ at $x$ lies in $Y_j$, and
write $x$ as $x_i^j$ if that neighboring vertex is $y_i^j$.  Since each vertex
of $Y$ has one neighbor via $M$, we have labeled $X$ so that
$X'=\{x_i^j\st 1\le i\le k\text{ and }1\le j\le 3\}$.

We construct a proper path-factor of $G$, dealing separately with each
component $C$ of $F^*$.  From $C$ we generate paths in $G$ that together
cover $V(C)$, the neighbors of $V(C)$ via $M$, and the vertices of $X_0$ whose
neighborhoods lie in $V(C)$.  The construction depends on whether the
restriction of $\bT$ to $C$ has an independent transversal or a spread
transversal.  For simplicity of notation, we describe the construction in the
case that $F^*$ is connected.  In the general case, $V(C)$ is the union of
$T_i$ for $i$ in some subset of $\{1,\ldots,k\}$, and the construction in
Case~1 or Case~2 covers all vertices in $T_i\cup\{x_i^0,x_i^1,x_i^2,x_i^3\}$
for each such index $i$.

{\bf Case 1:} {\it $Y_1$ is an independent transversal}.  We specify $k$ paths
of lengths 4, 6, or 8, each containing one vertex of $X_0$.  
Consider the paths $\la x_i^2,y_i^2,x_i^0,y_i^3,x_i^3\ra$ for $1\le i\le k$.
These paths are disjoint and cover $V(G)-(X_1\cup Y_1)$.  The $2k$ endpoints of
these paths form $X_2\cup X_3$.  Each vertex $y_i^1$ of $Y_1$ has one incident
edge $y_i^1x$ having color $2$ under $c$.  Since $Y_1$ is an independent
transversal, this neighbor $x$ lies in $X_2\cup X_3$, not in $X_1$.  Extend the
original path of length 4 ending at $x$ by adding $xy_i^1$ and $y_i^1x_i^1$.
Altogether there are $k$ such extensions to absorb $Y_1\cup X_1$.  Each of the
original paths extends by zero or two edges at each end, so we have the factor
using paths of the desired lengths.

{\bf Case 2:} {\it $Y_1$ is a spread transversal}.  Again each path contains one
vertex of $X_0$, but now we may also use length 2.  Specify an orientation of
each cycle in $F$, and delete the incoming edge to each vertex of $Y_1$.
Since $Y_1$ is a spread transversal, each cycle is cut, and what remains of $F$
consists of $k$ disjoint paths $P_1, \dots, P_k$, starting at $\VEC{y^1}1k$,
respectively, each with length at most 3.  By the definition
of $G$, each edge in $P_i$ expands to a path of length 2 in $G$ having edges
of colors 1 and 2 under $c$, yielding paths of even length (at most 6)
ending in $Y_2\cup Y_3$.  If some $P_i$ ends at $y\in Y_2\cup Y_3$, then the
neighbor of $y$ in $M$ is not covered by any of these paths.  Thus extending
each path $P_i$ by adding $x_i^0y_i^1$ at the beginning and the edge $yx$
of $M$ at the end yields $k$ disjoint paths of lengths in $\{2,4,6,8\}$
that cover $V(G)$.
\end{proof}

This method for finding proper factors is robust, since any proper
$3$-edge-coloring of $H$ and any indexing of its colors can be used.
Care is needed, since there exist 2-regular graphs $F$ and triple
systems $\bT$ where no mixed transversal exists, as shown in our final example.

\begin{example}
First we construct $F_1$ with no independent transversal.  Let $k_1$ be a 
multiple of 6, and let $F_1$ consist of ${k_1}/{2}$ cycles of length 4
and ${k_1}/{3}$ cycles of length 3.  Name the 4-cycles as
$[y_{2i-1}^1,y_{2i}^1,y_{2i-1}^2,y_{2i}^2]$ for $1 \leq i \leq k_1/2$.
Name the 3-cycles as $[y_{6i-3}^3,y_{6i-1}^3,y_{6i+1}^3]$ and
$[y_{6i-4}^3,y_{6i-2}^3,y_{6i}^3]$ for $1 \leq i \leq k_1/6$ (with
$y_{k_1+1}=y_1$).  An independent partial transversal has at most one vertex
in each cycle, and hence the largest independent partial transversal has at
most $k_1/2+k_1/3$ elements.

Next we construct $F_2$ with no spread transversal; for clarity, we use
vertices $z_i^j$ instead of $y_i^j$.  Let $k_2$ be a multiple of $2$, and let
$F_2$ consist of ${3k_2}/{2}$ cycles of length 2.  Name the 2-cycles as
$[z_i^1,z_{i+1}^2]$ for $1 \leq i \leq k_2$ (where $z_{k+1}=z_1$) and
$[z_{2i-1}^3,z_{2i}^3]$ for $1 \leq i \leq k_2/2$.  A transversal has
only $k_2$ elements, but a spread transversal must have an element in
each of the $3k_2/2$ cycles.

Both $F_1^*$ and $F_2^*$ are connected graphs.  To construct an example with no
mixed transversal, we start with disjoint copies of $F_1$
and $F_2$, with $k_1=12$ and $k_2=8$.  Exchange vertex $y_1^1$ in $F_1$
for $z_1^1$ in $F_2$.  This creates a new graph $F$ such that $F^*$ is
connected; hence a mixed transversal must be an independent transversal
or a spread transversal.  A transversal can cover at most 11 of the
2-cycles from $F_2$ and thus cannot be spread.  On the other hand, at most
11 vertices of $F_1$ can be chosen for an independent transversal, since
at most one vertex from each 3-cycle and 4-cycle can be selected.
We conclude that there is no mixed transversal.
\qed
\end{example}


\begin{thebibliography}{99}
\frenchspacing
%{\small
\renewcommand{\baselinestretch}{.92}            %line-spacing
\parskip=.5ex

\bibitem{AsCas}
A. S. Asratian and C. J. Cassegren,
On interval edge colorings of $(\alpha,\beta)$-biregular bipartite graphs,
{\it Discr.\ Math.{}} (in press).

%Some results on interval edge colorings of $(\alpha,\beta)$-biregular bipartite
%graphs, Research report LiTH-MAT-R-2006-09, Link\"oping University,
%Link\"oping, Sweden (2006).

\bibitem{AK1}
A. S. Asratian and R. R. Kamalian,
Interval coloring of the edges of a graph (in Russian).
{\it Applied Math.{}} 5 (1987), 25--34. 

\bibitem{AK2}
A. S. Asratian and R. R. Kamalian,
Investigation on interval edge-colorings of graphs.
{\it J. Combin.\ Theory Ser.\ B} 62 (1994), 34--43. 

\bibitem{Axe}
M. A. Axenovich,
On interval colorings of planar graphs.
Proc.\ 33rd Southeastern Intl.\ Conf.\  Combin.,
Graph Theory and Computing (Boca Raton, FL, 2002).
{\it Congr. Numer.{}} 159 (2002), 77--94. 

\bibitem{Cass}
C. J. Casselgren, Some results on interval edge colorings of bipartite graphs.
Master's Thesis, Link\"oping University, Link\"oping, Sweden (2005).

\bibitem{Giaro}
K. Giaro, The complexity of consecutive $\Delta$-coloring of bipartite graphs:
$4$ is easy, $5$ is hard.  {\it Ars Combin.{}}  47  (1997), 287--298.

\bibitem{GK1}
K. Giaro and M. Kubale,
Consecutive edge-colorings of complete and incomplete Cartesian products of
graphs.
Proc.\ 28th Southeastern Intl.\ Conf.\  Combin.,
Graph Theory and Computing (Boca Raton, FL, 1997).
{\it Congr. Numer.{}} 128 (1997), 143--149.

\bibitem{GK}
K. Giaro and M. Kubale, Compact scheduling of zero-one time operations in 
multi-stage systems, {\it Discrete Appl. Math.} 145 (2004), 95--103.

\bibitem{Hansen}
H. M. Hansen, Scheduling with minimum waiting periods (in Danish),
Master's Thesis, Odense University, Odense, Denmark (1992).

\bibitem{HansonLoten}
D. Hanson and C. O. M. Loten, 
A lower bound for interval colouring bi-regular bipartite graphs.
{\it Bull.\ Inst.\ Combin.\ Appl.{}}  18  (1996), 69--74.

\bibitem{HansonLotenToft}
D. Hanson, C. O. M. Loten, and B. Toft,
On interval colourings of bi-regular bipartite graphs.
{\it Ars Combin.{}} 50  (1998), 23--32.

\bibitem{Kamalian}
R. R. Kamalian,
Interval edge-colorings of graphs.
Doctoral Thesis, Novosibirsk (1990).

\bibitem{Kos}
A. V. Kostochka, unpublished manuscript, 1995.

\bibitem{Pyatkin}
A. V. Pyatkin, Interval coloring of $(3,4)$-biregular bipartite graphs having
large cubic subgraphs.  {\it J. Graph Theory}  47  (2004),  122--128.

\bibitem{Sevast}
S. V. Sevastjanov, 
Interval colorability of the edges of a bipartite graph.
{\it  Metody Diskret.\ Analiz.{}} 50 (1990), 61--72, 86.


\end{thebibliography}
\end{document}